\documentclass {amsart}
\usepackage[english]{babel}
\usepackage{amsthm}
\usepackage{amsmath}
\usepackage{amsfonts}

\def\R{\mathbb{R}}

\def\T{\mathbb{T}}
\def\N{\mathbb{N}}
\def\Z{\mathbb{Z}}
\def\cb{\mathcal{B}}

\def\cf{\mathcal{F}}

\def\ct{\mathcal{T}}

\def\eps{\varepsilon}

\newcommand{\Prob}{\mathbb P}
\newcommand{\mes}{\operatorname{mes}}
\renewcommand{\Im}{{\operatorname{Im}}}

\theoremstyle{definition}
\newtheorem{theorem}{Theorem}
\newtheorem{lemma}{Lemma}

\newtheorem{definition}{Definition}

\newtheorem{remark}{Remark}
\newtheorem{corollary}{Corollary}

\newcommand{\ld}[1]{\operatorname{LD}(#1)}
\newcommand{\bd}[1]{\operatorname{BD}(#1)}
\newcommand{\od}[1]{\operatorname{OD}(#1)}
\renewcommand{\sc}[1]{\operatorname{SC}(#1)}

\newcommand{\interior}{\operatorname{Int}}
\newcommand{\closure}{\operatorname{Cl}}

\newcommand{\xor}{\triangle}

\author{Ron Peled}
\title{Restoring Topology from Shifts}
\begin{document}
\begin{abstract}
It is known that the topology of a Polish group is uniquely determined by its Borel structure and group operations, but this does not give us a way to find the topology. In this article we expand on this theorem and give a criterion for a measurable function on the group to be continuous. We show that it is continuous iff there exists some second countable topology in which all the shifts of the function are continuous. Here measurability can be taken to be measurable with respect to Haar measure (in locally compact groups), or having the Baire property or universally measurable (in Abelian Polish groups). Our results appear to be new even when the group is $\R$. As a special case we get that a measurable homomorphism is continuous (a known result).

As an application, We give a first proof of a dichotomy guessed by Tsirelson long ago for stationary stochastic processes. Either the process is sample continuous, or its paths cannot be continuous in any second countable topology on any non-null set.
\end{abstract}
\maketitle
\begin{section}{Introduction}
It is well known  that the topology of a Polish group is uniquely determined by its Borel structure and group operations (see e.g. proposition 12.24 in \cite{KE1}), but this does not provide us with a way of finding this topology.
In this article we extend this result by giving a criterion for a measurable function on the group to be continuous by examining only the shifts of the function as functions on an abstract set. Our criterion is another manifestation of the well observed principle that if a function has some algebraic properties and some analytic properties then it must already be continuous (for example, a measurable group homomorphism is continuous, see theorem 9.10 in \cite{KE1}, in many situations, this example is a special case of our theorem), our criterion asserts that if there exists some second countable topology on the group in which all the shifts of the function are continuous, then the function is actually continuous in the usual metric of the group. Our results appear to be new even when the group is $\R$.

This criterion also has a local version stating that if the restriction of many (in some sense) small shifts of the function to a subspace containing a large portion (in some sense) of the neighborhood of a point are continuous in the relative second countable topology then the function is continuous in a neighborhood of that point (theorems \ref{local_measurable_thm} and \ref{local_category_thm}). This can even be weakened further by allowing each shift of the function to be changed on a null set and requiring only that these changed shifts are continuous in the second countable topology, even so we conclude that the function is almost everywhere equal to a continuous function (theorems \ref{local_measurable_process_thm} and \ref{local_category_process_thm}).

As in similar results, our criterion turns out to be valid in more than one sense of measurability, the conclusion holds if our function is measurable with respect to Haar measure (on a locally compact Polish group, theorems \ref{general_local_measurable_thm} and \ref{general_local_measurable_process_thm}) or when it has the Baire property (on any Polish group, theorems \ref{general_local_category_thm} and \ref{local_category_process_polish_groups_thm}), or when it is universally measurable (on any Abelian Polish group, unfortunately we do not know how to handle the non-Abelian case, theorem \ref{global_universally_measurable_thm}). For the local version of our criterion we require the notion of a density point of a set (or function), such a notion is known for functions on $\R$ or for functions with the Baire property. We additionally give another such notion for a function on a general locally compact metric group.

This work began after the proof of an old claim by Tsirelson (Corollary 1 in \cite{TS1}) was found to be false, and in fact, many of the ideas in this work are based on ideas from that article. We also give a correct proof of that claim which gives a dichotomy for stationary stochastic fields on locally compact Polish time sets (this also appears to be new even for $\R$). We prove that either the process is sample continuous or its paths cannot be continuous in any second countable topology on any non-null set (theorems \ref{stationary_process_thm} and \ref{stationary_process_locally_compact_thm}). It is of interest to compare this with Belyaev's dichotomy \cite{BE1} for \emph{Gaussian} stationary processes which states that such a process is either sample continuous, or almost all of its paths are unbounded in every interval.

The author is very grateful to Prof. B. S. Tsirelson for asking the initial question and providing many helpful ideas, insights and remarks without which this work would never have been completed. 

The author would also like to thank Prof. A. Olevskii who gave the starting point for this research \cite{OL1} by answering a question asked by Tsirelson and thus providing the main tool in the first proof of theorem \ref{local_measurable_thm}. The theorem was later found to have a simpler proof which is given in the article, but the original theorem of Olevskii is still given here in an appendix (theorem \ref{Olevskii_theorem}) with his kind permission since the technique used has its own independent interest.
\end{section}
\begin{section}{Definitions, theorems and remarks} \label{real_formulation_section}
In the following theorems, the definitions are given in full generality, but for simplification the formulation of the theorems and the proofs are done in $\R$. In section \ref{general_formulation_section} we will give the more general formulation of the theorems for Polish groups and indicate what changes (if any) are needed for the proofs. Those results which use Lebesgue density points are specific to $\R$, but we will also give in section \ref{general_formulation_section} a similar notion for a general locally compact Polish group and get similar results there.

In order to state our theorems we need a definition:
\begin{definition}
Let $G$ be a group, $H$ a topological space and $f:G\to H$ a function. Given sets $S,T\subset G$, we call $f$ \emph{$(S,T)$ shift continuous} if there exists a second countable topology $\ct$ on $T$ such that for every $s_0\in S$, $f(s_0+\cdot)|_{T}$ is $\ct$-continuous.
\end{definition}
It is easy to see that if $G$ admits a separable metrizable topology making it a topological group then existence of a topology $\ct$ as in the definition implies existence of a separable metric topology $\ct'$ also good for the definition\footnote{To see this, take $\{U_n\}_n$ a basis for $\ct$, consider $\alpha:G\to\{0,1\}^\N$, $\alpha(x):=(1_{U_1}(x),1_{U_2}(x),\ldots)$, transport the metric of $\{0,1\}^\N$ through $\alpha$ to $G$. This gives a second countable semi-metric topology. Add this semi-metric to the metric of $G$ to make it into a separable metrizable topology.}. In our theorems $G$ will always be a Polish group.

Note that if $G$ is not Abelian, this should be called left shift continuous and an analogous concept of right shift continuous should be defined. Our theorems are true for both definitions with the same proofs and we consider only the definition as given above.

An equivalent definition which does not involve abstract topologies on the group is\footnote{We will not use this definition in our theorems or proofs.}:
\begin{lemma} \label{equiv_def_lemma}
Let $G$ be a second countable topological group, $H$ a topological space, for a function $f:G\to H$:
\begin{enumerate}
\item If there exists a second countable space $X$ and a function $h:X\to G$ with $T\subset\Im(h)$ such that for every $s_0\in S$, $f(s_0+h(\cdot))$ is continuous then $f$ is $(S,T)$ shift continuous.
\item if $f$ is $(S,T)$ shift continuous then there exists a second countable space $X$ and a continuous one-to-one function $h:X\to G$ with $T=\Im(h)$ such that for every $s_0\in S$, $f(s_0+h(\cdot))$ is continuous.
\end{enumerate}
\end{lemma}
In what follows, we will examine the value of a function $f:\R\to\R$ at the point $s+t$. In order to understand the theorems better, it may be helpful for the reader to have in mind the two variable function $g(t,s):=f(t+s)$ defined on $\R^2$ with the $t$ axis and the $s$ axis. 

We will consider sets $S,T\subset\R$ which are big in some sense and in this case we will show that shift continuity is equivalent to continuity in the standard metric. It turns out that our theorems are true when the sets are big either in Lebesgue measure sense (more generally, Haar measure sense) or in Baire category sense. 

We write $\mes$ and $\mes^*$ for Lebesgue measure and outer Lebesgue measure on $\R$ respectively, we will use the notion of Lebesgue density point\footnote{We will sometimes abbreviate this term to Lebesgue point.} and outer density point (that is, a Lebesgue density point with respect to outer measure $\mes^*$), for a given set $S\subset\R$ we define $\ld{S}$ to be its set of Lebesgue density points and $\od{S}$ its set of outer density points (some may be outside of $S$). We remind the reader that by the outer density and Lebesgue density theorems (thm 3.20 p. 17 in \cite{MC1}), for any set $S$, $S\setminus \od{S}$ is a null set.\footnote{In this article, we may replace $\od{S}$ with the slightly larger set $\ld{S^c}^c$ and the proofs remain the same.} and if $S$ is measurable then $S\xor\ld{S}$ is a null set. For ease of reference, we collect all the notations here (Baire category definitions are given below):

\framebox[\width][c]{
\begin{tabular}{l}
$\mes$ - Lebesgue measure\\
$\mes^*$ - outer Lebesgue measure\\
$\ld{S}$ - Lebesgue density points of $S$\\
$\od{S}$ - outer Lebesgue density points of $S$\\
$\bd{S}$ - Baire density points of $S$\\
$\sc{S}$ - second category points of $S$
\end{tabular}}

For sets $A,B\subset\R$ we write $A+B$ for their Minkowski sum $\{a+b\ |\ a\in A, b\in B\}$. We will prove:
\begin{theorem}\label{local_measurable_thm}
Let $f:\R\to\R$ be a measurable function. If $f$ is $(S,T)$ shift continuous for a measurable $S\subset\R$ and a set $T\subset\R$ then $f$ is continuous at all points of $\ld{S}+\od{T}$\footnote{Note that this is an open set from Steinhaus' lemma (\cite{ST1} or thm 4.8 p. 21 and p. 93 of \cite{MC1}).}.
\end{theorem}

For the Baire category version of the theorem, we first define:
\begin{definition}
Let $X$ be a topological space, $S\subset X$ and $s\in X$. We call $s$ a \emph{Baire density point} of $S$ if there exists an open set $U$ containing $s$ such that $U\setminus S$ is of first category.
\end{definition}
\begin{definition}
Let $X$ be a topological space, $S\subset X$ and $s\in X$. We call $s$ a \emph{second category point} of $S$ if for every open set $U$ containing $s$, $S\cap U$ is of second category (or equivalently, if $s$ is not a Baire density point of $S^c$).
\end{definition}
For given $S\subset X$ we define $\bd{S}$ to be its set of Baire density points and $\sc{S}$ its set of second category points (some may be outside of $S$).  Clearly $\bd{S}\subset\sc{S}$, we note that $\bd{S}$ is open and remind the reader that $S$ has the Baire property iff $\bd{S}\xor S$ is a first category set (see, for example, \cite{MC1}). Note also that from the Banach category theorem (thm 16.1 p. 62 and p. 96-97 of \cite{MC1}), for any $S$, $S\setminus \sc{S}$ is of the first category (see, for example, thm 8.29 p. 49 of \cite{KE1} or lemma 8.4 p. 100 of \cite{HO1}\footnote{\label{heavy_point_footnote}In \cite{HO1} p. 100 they define the set of \emph{heavy} points of $S$ to be $\interior(\sc{S})$. This is in general a set larger than $\bd{S}$ (consider a Bernstein set where every point is a heavy point, but no point is a Baire density point), but it is easy to check that they coincide for sets $S$ with the Baire property.}).

We will call a set $S\subset X$ \emph{residual} if $X\setminus S$ is of first category.

The Baire category version of theorem \ref{local_measurable_thm} is:
\begin{theorem}\label{local_category_thm}
Let $f:\R\to\R$ be a function with the Baire property. If $f$ is $(S,T)$ shift continuous for $S\subset\R$ with the Baire property and a set $T\subset\R$ then $f$ is continuous at all points of $\bd{S}+\sc{T}$.\footnote{Which is an open set since $\bd{S}$ is open.}
\end{theorem}

\begin{remark} \label{large_category_remark}
We note that here and in theorem \ref{local_measurable_thm} if $S$ and $T$ are big enough, the conclusion is that $f$ is continuous on the whole of $\R$. For example (in theorem \ref{local_category_thm}), if $S$ is residual and $T$ is of second category (because this implies that $\sc{T}$ is of second category) or $S$ is of second category (and has the Baire property) and $T$ is residual (or even only $\sc{T}=\R$).
\end{remark}

We also give a theorem similar to theorem \ref{local_measurable_thm} in the language of stochastic processes which we phrase as a dichotomy\footnote{The reader not interested in stochastic processes should jump to remark \ref{measurable_process_remark} and then to theorem \ref{local_category_process_thm}.}. We need the definition of \emph{natural modification} from \cite{TS1}, this is the main definition from that paper and it is very similar to our shift continuity. We remind the reader of the settings there, we consider a stochastic process $(X_t)_{t\in T}$ on the standard Lebesgue probability space $(\Omega,\cf,\Prob)$ and the time set $T$.
\begin{definition}
A modification $\xi:\Omega\to \R^T$ is called \emph{natural} if there exists a second countable topology $\ct$ on $T$ such that for a.e. $\omega\in\Omega$, $\xi(\omega)$ is $\ct$-continuous.
\end{definition}
The definition given in \cite{TS1} is different in that it requires $\ct$ to be separable metrizable, but it is easily seen that all the results there remain true when using this apparently weaker definition instead (and in particular, existence of a second countable topology $\ct$ as in our definition implies existence of a separable metrizable topology $\ct'$ as in \cite{TS1}, hence the two definitions are equivalent).

The reader should note the natural modification is in many senses a canonical modification of the process (hence the name natural), for many equivalent definitions and properties see theorems 1 and 2 in \cite{TS1}.

We will denote by $(X_t)|_A$ the restriction of the stochastic process $(X_t)_{t\in T}$ to a smaller time set $A\subset T$.
\begin{theorem}  \label{local_measurable_process_thm}
Let $f:\R\to\R$ be a measurable function, $\Prob$ a probability measure on $\R$ equivalent to Lebesgue measure. Define a stochastic process $(X_t)_{t\in\R}$ on $(\R,\cf,\Prob)$ by $X_t(\omega):=f(\omega+t)$. Then there is an alternative:

Either $(X_t)_{t\in\R}$ is sample continuous or for every non-null set $T\subset\R$, $(X_t)|_T$ does not have a natural modification.
\end{theorem}
\begin{remark} \label{sample_continuity_remark}
Note that $(X_t)_{t\in\R}$ is sample continuous iff $f$ is a.e. equal to a continuous function.

Also, if $(X_t)_{t\in\R}$ is sample continuous and for a set $T\in\R$, $(X_t)|_T$ has a natural modification then, from the uniqueness of the natural modification (theorem 2(a) in \cite{TS1}), for a.e. $\omega\in\Omega$ the sample continuous modification and the natural modification coincide on $T$.
\end{remark}
\begin{remark} \label{measurable_process_remark}
We will actually prove a generalization of this which does not use all the shifts of the function. This version is easier to formulate without the language of stochastic processes and not in the form of an alternative:

Let $f:\R\to\R$ be a measurable function. Given a measurable $S\subset\R$, a set $T\subset\R$ of positive outer measure and a function $g:(S\times T)\to\R$ such that for each $t\in T$, $g(\cdot,t)$ is a.e. equal to $f(\cdot+t)$. If there exists a second countable topology $\ct$ on $T$ such that for each $s\in S$, $g(s,\cdot)$ is $\ct$-continuous, then $f$ is a.e. equal to a function $h$ continuous on the open set $\ld{S}+\od{T}$ and for a.e. $s\in S$, $g(s,\cdot)=h(s+\cdot)$.

It is not hard to check that theorem \ref{local_measurable_thm} is the special case where for each $t\in T$, $g(\cdot,t)=f(\cdot+t)$ (and not just a.e. equal).
\end{remark}

We continue to give a corrected proof of the corollary from Tsirelson's paper (corollary 1 in \cite{TS1}) and extend the result somewhat. The following definition doesn't appear to be in wide use, so we give it here for clarity (see also p. 420 \cite{DU1}):
\begin{definition}
Given a process $(X_t)_{t\in T}$ where the time set $T$ is a measurable space with $\sigma$-algebra $\cb$. We call $(X_t)_{t\in T}$ a \emph{measurable} process if it has a modification $f(\omega,t)$ measurable in both variables jointly with respect to $\cf\otimes\cb$.
\end{definition}
\begin{remark} \label{evaluation_function_remark}
There is another way to define a measurable process without using modifications. Consider the space of equivalence classes (mod 0) of random variables on $\Omega$ equipped with the (metrizable) topology of convergence in probability, we denote this space by $L_0(\Omega)$ (see also p. 226 \cite{DU1}). One can also think of a process (without specifying a modification of it) as a mapping from $T$ to $L_0(\Omega)$. It turns out that this mapping is measurable if and only if the process is measurable.

One direction of the equivalence is obvious. The other requires the notion of an "evaluation function", a measurable function $\alpha:\Omega\times L_0(\Omega)$ which satisfies $\alpha(\cdot,f)=f$ for each $f\in L_0(\Omega)$\footnote{Many such functions exist, for example (taken from intro. of \cite{GTW1}) $\alpha(\omega,f):=\tan(\limsup_{\eps\to0}\frac{1}{2\eps}\int_{\omega-\eps}^{\omega+\eps}\arctan(f)(x_1)dx_1)$.}. Given such a function we immediately obtain a measurable modification for the process from a measurable mapping of $T$ to $L_0(\Omega)$.
\end{remark}

\begin{theorem}\label{stationary_process_thm}
Let $(X_t)_{t\in\R}$ be a measurable \emph{stationary} stochastic process. Then the following alternative holds:

Either $(X_t)_{t\in\R}$ is sample continuous or for every non-null set $T\subset\R$, $(X_t)|_T$ does not have a natural modification.
\end{theorem}
\begin{remark} \label{second_sample_continuity_remark}
As in theorem \ref{local_measurable_process_thm}, if $(X_t)_{t\in\R}$ is sample continuous and for a set $T\in\R$, $(X_t)|_T$ has a natural modification then, from the uniqueness of the natural modification (theorem 2(a) in \cite{TS1}), for a.e. $\omega\in\Omega$ the sample continuous modification and the natural modification coincide on $T$.
\end{remark}
\begin{remark}
It is known that a measurable stationary process on $\R$ is automatically continuous in probability (that is, the mapping from $T$ to $L_0(\Omega)$ is continuous). We will elaborate more on this in the proof of the theorem.
\end{remark}

Finally, we end by proving a Baire category version of theorem \ref{local_measurable_process_thm}:
\begin{theorem}  \label{local_category_process_thm}
Let $f:\R\to\R$ be a function with the Baire property. Given a set $T\subset\R$ of second category and a function $g:\R^2\to\R$ such that for each $t\in T$ $g(\cdot,t)$ is equal to $f(\cdot+t)$ on a residual set, then if  there exists a second countable topology $\ct$ on $T$ such that for a residual set of $s\in\R$, $g(s,\cdot)|_T$ is $\ct$-continuous, then $f$ is equal on a residual set to a continuous function $h$ and for a residual set of $s\in \R$, $g(s,\cdot)|_T=h(s+\cdot)|_T$.
\end{theorem}
\begin{remark} \label{local_category_process_remark}
A similar generalization to that in remark \ref{measurable_process_remark} is also true here (it is formulated precisely in the proof section).
\end{remark}
\end{section}

\begin{section}{Proofs and some remarks}
We begin by proving lemma \ref{equiv_def_lemma}:
\begin{proof} (of \underline{lemma \ref{equiv_def_lemma}})
\begin{enumerate}
\item Assume that there exists a second countable space $X$ and a function $h:X\to G$ such that for each $s_0\in S$, $f(s_0+h(\cdot))$ is continuous, we will show that $f$ is $(S,\Im(h))$ shift continuous.

Let $\{U_n\}$ be a countable basis for the topology of $X$ and define the topology $\ct$ on $\Im(h)$ by setting $V_n:=h(U_n)$ as its basis. Using this definition $h$ is an open mapping, hence it is clear that $f$ is $(S,\Im(h))$ shift continuous with the $\ct$ topology. Of course this implies that $f$ is also $(S,T)$ shift continuous for any $T\subset\Im(h)$.
\item If $f$ is $(S,T)$ shift continuous for some $S,T\subset G$, take $X=T$ with the second countable topology $\ct$ of the shift continuity, and the function $h(x)=x$. To make $h$ continuous, enrich the topology on $X$ to the weakest topology which also contains that of the subspace $T$ of $G$. This leaves $X$ second countable since $G$ is second countable.
\end{enumerate}
\end{proof}

\begin{subsection}{Baire category case}

The Baire category case is simpler than the measurable case and illustrates some of the main ideas, we begin with the proof of theorem \ref{local_category_process_thm} and get theorem \ref{local_category_thm} as a conclusion.

We begin with a definition and a lemma:
\begin{definition}
Let $X$ be metric space, $s\in X$ is a \emph{Baire density point} of a function $f:X\to\R$ if for each open $U\subset \R$ with $f(s)\in U$, $s$ is a Baire density point of $f^{-1}(U)$.
\end{definition}
We note that from theorem 8.1 in \cite{MC1}, if $f$ has the Baire property then it has a residual set of Baire density points.\footnote{It is even true (lemma 8.5 p. 100 of \cite{HO1}) that  \emph{every} $f:X\to\R$ has a residual set of heavy points (or almost continuity points), where $s$ is a heavy point of $f$ if for every open $U$ with $f(s)\in U$, $s$ is a heavy point (footnote \ref{heavy_point_footnote}) of $f^{-1}(U)$.}.

We will now show how to choose a special representative for the class of changes of a function $f$ (with the Baire property) on a first category set. For each $f$ with the Baire property, let $\bd{f}$ be its (residual) set of Baire density points and define a function $\tilde{f}$ by:
\[\tilde{f}(x):=\sup\{\limsup f(x_n)\ |\ \{x_n\}_n\subset \bd{f}, \lim x_n=x\}\]
Where we mean that $\tilde{f}$ can assume the values $\pm\infty$. It is easily seen from the definition that $f=\tilde{f}$ on the set $\bd{f}$ since $f|_{\bd{f}}$ is continuous. It follows that $\tilde{f}$ differs from $f$ on a set of first category (and in particular, $\tilde{f}$ is finite on a residual set and is Baire measurable).

It is also easily checked that if $g$ is a change of $f$ on a first category set then $\tilde{f}=\tilde{g}$ (hence the term representative).

The following lemma is the main property of the special representative that we will use (it should be noted, however, that there is more than one way to define a representative with this property):
\begin{lemma} \label{Baire_canonical_representative_lemma}
If there exists a set $S\subset\R$ and a point $x\in (\bd{S}\cap S)$ such that $f|_S$ is continuous at $x$, then $\tilde{f}$ is continuous at $x$ and $f(x)=\tilde{f}(x)$.
\end{lemma}
\begin{proof}
To prove continuity of $\tilde{f}$ at $x$, take any sequence $\{x_n\}_n$ with $x_n\to x$, we will show that $\lim \tilde{f}(x_n)=\tilde{f}(x)=f(x)$ (and in particular the limit exists). For each $n$, from the definition of $\tilde{f}(x_n)$ we can take $y_n\in \bd{f}$ such that $|y_n-x_n|<\frac{1}{n}$ and $|\tilde{f}(x_n)-f(y_n)|<\frac{1}{n}$ (if $\tilde{f}(x_n)=\pm\infty$, take $y_n$ with $|f(y_n)|>n$). It is enough to show that $\lim f(y_n)=f(x)$ (since $f(x)$ is finite, this implies in particular that only for a finite number of $x_n$, $\tilde{f}(x_n)=\pm\infty$). Now, since $y_n\in \bd{f}$, $y_n\to x$ and $x\in(\bd{S}\cap S)$ we can choose a sequence $\{z_n\}\subset (\bd{f}\cap(\bd{S}\cap S))$ such that for large enough $n$, each $z_n$ is a small perturbation of $y_n$, still $z_n\to x$ and $\lim f(z_n)= \lim f(y_n)$ (in the sense that they exist together and are equal). Now, from the continuity of $f|_S$ we get $\lim f(z_n)=f(x)$ and we are done.
\end{proof}

Following Tsirelson \cite{TS1} (proof of theorem 1 part (a)$\to$uniform Lusin measurability) for a similar setting in a measure space, we continue by proving:
\begin{lemma} \label{uniform_lusin_category_lemma}
Given $S\subset\R$ with the Baire property, $T\subset\R$ a set and $g:S\times T\to\R$ such that for each $t\in T$, $g(\cdot, t)$ has the Baire property, then if  there exists a second countable topology $\ct$ on $T$ such that for each $s\in\R$, $g(s,\cdot)$ is $\ct$-continuous then there exists a set $S_1\subset S$ with $S\setminus S_1$ of first category (in $\R$) such that for each $t\in T$, $g(\cdot,t)|_{S_1}$ is continuous (in the standard metric of $\R$).
\end{lemma}
\begin{proof} Let $\{U_n\}$ be a countable basis for the $\ct$ topology on $T$, define for each $n\in\mathbb{N}$ the functions $g_n^{\sup}(s):=\sup_{t\in U_n}g(s,t)$ and $g_n^{\inf}(s):=\inf_{t\in U_n}g(s,t)$ then $g_n^{\sup}$ and $g_n^{\inf}$ have the Baire property since it is enough to take the $\sup$ and $\inf$ on a $\rho$-dense set, it follows from theorem 8.1 (p. 36) in \cite{MC1} (applied to a countable set of functions) that there exists a residual set $S_1\subset S$ with $S\setminus S_1$ of first category\footnote{\label{first_category_footnote}Note that since $S$ has the Baire property, if $M\subset S$ is of first category in $S$, then it is of first category in $\R$ because $S$ can be written as an $F_\sigma$ set minus a first category set (thm 4.4 p. 21 in \cite{MC1}).} such that all of the $g_n^{\sup}$ and $g_n^{\inf}$ restricted to $S_1$ are continuous. We claim that this is the required set.

To see this, fix $t_1\in T$ and let us show continuity of $g(\cdot,t_1)|_{S_1}$. Fix $s_1\in S_1$ and $\epsilon>0$. Let $n$ be such that $t_1\in U_n$ and:
\[\forall t\in U_n \qquad |g(s_1,t)-g(s_1,t_1)|<\epsilon\]
then $g_n^{\sup}(s_1)-\epsilon\le g(s_1,t_1)\le g_n^{\inf}(s_1)+\epsilon$. Now from continuity of $g_n^{\sup}(s)$ and $g_n^{\inf}(s)$ restricted to $S_1$, it follows that for $s\in S_1$ close enough to $s_1$ we have 
\[g(s_1,t_1)-2\epsilon<g_n^{\inf}(s)\le g(s,t_1)\le g_n^{\sup}(s)<g(s_1,t_1)+2\epsilon\]
\end{proof}
We now proceed to the proof of the theorem:
\begin{proof} (of \underline{theorem \ref{local_category_process_thm}}) Applying lemma \ref{uniform_lusin_category_lemma} to the function $g$, we get a residual set $S_1\subset\R$ with the property of the lemma.

The required function $h$ is simply $\tilde{f}$ (and $\tilde{f}$ is finite everywhere). To see this, fix $t\in T$, since $g(\cdot,t)$ is a change of $f(\cdot+t)$ on a first category set and $g(\cdot,t)|_{S_1}$ is continuous we get from lemma \ref{Baire_canonical_representative_lemma} (since $\bd{S_1}\cap S_1=S_1$) that $h(\cdot+t)=g(\cdot,t)$ on $S_1$ and $h(\cdot+t)$ is continuous on $S_1$. We conclude that $h$ is continuous at $S_1+T$ and since $T$ is of second category $S_1+T=\R$ and $h$ is continuous everywhere as required.
\end{proof}
Looking carefully at the proof of theorem \ref{local_category_process_thm} we see that actually the following generalization is true:
\begin{theorem}  \label{generalized_local_category_process_thm}
Let $f:\R\to\R$ be a function with the Baire property. Given a set $S\subset\R$ with the Baire property, a set $T\subset\R$ and a function $g:(S\times T)\to\R$ such that for each $t\in T$, $g(\cdot,t)$ is equal on a residual set\footnote{See footnote \ref{first_category_footnote}.} to $f(\cdot+t)$. If  there exists a second countable topology $\ct$ on $T$ such that for each $s\in\R$, $g(s,\cdot)$ is $\ct$-continuous, then $f$ is equal on a residual set  to a function $h$ continuous on $S_0+T$ where $S_0\subset S$ and $S\setminus S_0$ is of first category. Furthermore, $g(s,\cdot)=h(s+\cdot)$ for each $s\in S_0$.
\end{theorem}
The proof is the same as that of theorem \ref{local_category_process_thm}. The set $S_0$ in its conclusion is the set $S_1\cap\bd{S_1}$ (where $S_1$ is the set obtained from lemma \ref{uniform_lusin_category_lemma}).

We deduce theorem \ref{local_category_thm}:
\begin{proof} (of \underline{theorem \ref{local_category_thm}}) 
Define $g:(S\times T)\to\R$ by $g(s,t):=f(s+t)$ then the conditions of theorem \ref{generalized_local_category_process_thm} are fulfilled. It follows that there exists a function $h$ continuous on $S_0+T$ for $S\setminus S_0$ of first category and such that $f|_{S_0+T}=h|_{S_0+T}$ (using the furthermore part of the theorem). Now noting that $S_0+T\supseteq \bd{S}+\sc{T}$ and that $\bd{S}+\sc{T}$ is open (since $\bd{S}$ is open) we get that $f$ is continuous at $\bd{S}+\sc{T}$ as we sought to prove.
\end{proof}

\begin{remark}
We allowed $T$ in Theorem \ref{generalized_local_category_process_thm} to be of first category. Note that in this case, not much information is gained about $f$ since we know that for any function $f$ with the Baire property, the function $\tilde{f}$ is continuous on a residual set, namely that of its Baire density points. Usually, the interesting case is only when the continuity set of $\tilde{f}$ contains interior points (for example, only then can we deduce theorem \ref{local_category_thm}). However, there is some interest in the furthermore part of the theorem which shows that the sets $S\setminus S_t$ on which $g$ differs from $f$ are actually all contained in a single first category set.
\end{remark}

\end{subsection}
\begin{subsection}{Measurable case}

We now proceed to the measurable case, i.e., to (Lebesgue) measurable shift continuous functions.

We start with a lemma analogous to lemma \ref{uniform_lusin_category_lemma}:
\begin{lemma} \label{uniform_lusin_measurable_lemma}
Given $S\subset\R$ a measurable set, $T\subset\R$ a set, $\mu$ a probability measure on $S$ and $g:S\times T\to\R$ such that for each $t\in T$, $g(\cdot,t)$ is $\mu$-measurable\footnote{Measurable for the completion of $\mu$.}, then if there exists a second countable topology $\ct$ on $T$ such that for each $s\in S$, $g(s,\cdot)$ is $\ct$-continuous then the functions $\{g(\cdot,t)\}_{t\in T}$ are \emph{uniformly Lusin measurable} in the sense that for each $\eps>0$ there exists a compact $K_{\eps}\subset S$ with $\mu(S\setminus K_{\eps})<\eps$ such that for all $t\in T$, $g(\cdot,t)|_{K_\eps}$ is continuous.
\end{lemma}
\begin{proof}
Let $\mu$ be a measure as in the lemma, aiming to use Tsirelson's results in \cite{TS1}, we define a stochastic process $X_t$ on $(S,\mu)$ and the time set $T$ by $X_t(s):=g(s,t)$. From the assumptions of the lemma, we see that $g$ is a natural modification for $X_t$, hence as was proven in \cite{TS1} (proof of theorem 1 part (a)$\to$uniform Lusin measurability) we have that for every $\eps>0$ there exists a compact $K_{\eps}\subset S$ with $\mu(S\setminus K_{\eps})<\eps$ such that for all $t\in T$, $g(\cdot,t)|_{K_\eps}$ is continuous (in \cite{TS1} this is shown when the probability space is $([0,1],\mes)$, but the proof goes through whenever the probability space is a measurable subset of a Polish space since only Lusin's theorem (thm 17.12 p. 108, ex. 17.15 p. 109 and thm 13.1 p. 82 in \cite{KE1}) is required, see also lemma \ref{general_uniform_lusin_measurable_lemma} for a proof in more general settings).
\end{proof}

We would like to proceed in the same way as in the Baire category case by defining a special representative (for the class of measure zero changes of a measurable function). We can still use the same definition for that representative (replacing Baire density points with Lebesgue density points), but unfortunately, for such a representative the analogue of lemma \ref{Baire_canonical_representative_lemma} for Lebesgue density points is no longer true.  Instead, in the following proof, the main tool we will use is that given a measurable set $A$ with finite measure, the function $\mes((x+A)\cap A)$ is continuous in $x$.

We continue to prove theorem \ref{local_measurable_thm}:
\begin{proof} (of \underline{theorem \ref{local_measurable_thm}}) Define $G:S\times T\to\R$ by $g(s,t):=f(s+t)$, then the conditions of lemma \ref{uniform_lusin_measurable_lemma} are fulfilled. If $\mes(S)=0$ there is nothing to prove, otherwise fix a probability measure $\mu$ on $S$ equivalent to Lebesgue measure, it follows that for every $\eps>0$ there exists $K_\eps\subset S$ with $\mu(S\setminus K_\eps)<\eps$ (we don't use compactness of $K_\eps$) such that for all $t\in T$, $f(\cdot+t)|_{K_\eps}$ is continuous.

Now, take any point $x_0=s_0+t_0$ where $s_0\in \ld{S}$ and $t_0\in\od{T}$ and assume in order to get a contradiction that there exists a sequence $x_n\to x_0$ with $\lim f(x_n)\ne f(x_0)$. We note that it is sufficient to find $\delta>0$ and a $t_1\in T$ such that $(x_0-t_1)\in K_\delta$ and a subsequence $(x_{n_k}-t_1)\in K_\delta$ since this will contradict continuity of $f(\cdot+t_1)|_{K_\delta}$. Defining now for any $\delta>0$ and subsequence $(x_{n_k})_k$ the set $K_\delta^{n_k}:=(x_0-K_\delta)\cap(\cap_k(x_{n_k}-K_\delta))$ we need only find such a set with $T\cap K_\delta^{n_k}\neq\phi$.

Let $\eps>0$ (to be fixed later), denoting for any $a\in\R$, $I_\eps^a:=(a-\frac{\eps}{2},a+\frac{\eps}{2})$ we take $\delta>0$ so small that \[\mes(K_\delta \cap I_\eps^{s_0})\ge\frac{1}{2}\mes(S\cap I_\eps^{s_0})\]
We use the above-mentioned fact that for any measurable $A\subset \R$, the function $x\to\mes((x+A)\cap A)$ is continuous (if $\mes{A}=\infty$ we use that for $x$ small enough $\mes((x+A)\cap A) \ne 0$) and deduce that we can choose a subsequence $x_{n_k}$ such that \[\mes(K_\delta^{n_k}\cap I_\eps^{t_0})\ge\frac{1}{3}\mes(S\cap I_\eps^{s_0})\]
We claim that if $\eps$ is sufficiently small, $T\cap K_\delta^{n_k}$ is not empty as required. To see this, use that $s_0\in \ld{S}$ and $t_0\in \od{T}$ to choose $\eps$ such that 
\[\mes(S\cap I_\eps^{s_0})>\frac{3}{4}\eps \qquad \rm{ and } \qquad \mes^*(T\cap I_\eps^{t_0})>\frac{3}{4}\eps\]
Then $\mes(K_\delta^{n_k}\cap I_\eps^{t_0})>\frac{1}{4}\eps$ so that $T\cap K_\delta^{n_k}\cap I_\eps^{t_0}$ is not empty as required.
\end{proof}

We will generalize this proof to prove theorem \ref{local_measurable_process_thm}.  we first need a well known lemma\footnote{But since we couldn't find a reference to the proof, we give one here.}:
\begin{lemma} \label{ae_measurable_continuity_criterion}
Given a measurable $D\subset\R$, let $f: D\to\R$ be a locally bounded measurable function, if $f$ is not a.e. equal (on $D$) to a continuous function, then there exists a point $x_0\in D$ and sequences $x_n^1\to x_0$ and $x_n^2\to x_0$ of Lebesgue points of $f$ with $\lim f(x_n^1)\ne \lim f(x_n^2)$.
\end{lemma}
\begin{proof}
Assume the contrary, then since $f$ is locally bounded we get that for any $x_0\in D$ and any sequence $x_n\to x_0$ of Lebesgue points of $f$, the limit $f(x_n)$ exists and is independent of the choice of $(x_n)_n$. Hence, denoting the set of Lebesgue points of $f$ by $\ld{f}$, we can define as in the Baire category case, the function:
\[\tilde{f}(x):=\sup\{\limsup f(x_n)\ |\ \{x_n\}_n\subset \ld{f}, \lim x_n=x\}\]
and we get from the above that $\tilde{f}=f$ on $\ld{f}$ (so it is a change of $f$ on a measure zero set) and that for any $x_0\in D$, $\tilde{f}|_{\ld{f}\cup{x_0}}$ is continuous. It is now easy to conclude that $\tilde{f}$ is actually continuous on $D$ since for any $x_0\in D$ and $x_n\to x_0$, we can find (from the definition of $\tilde{f}$) a sequence $\{y_n\}_n\subset \ld{f}$ with $y_n\to x_0$ and $\lim \tilde{f}(y_n)=\lim \tilde{f}(x_n)$ (like in the proof of lemma \ref{Baire_canonical_representative_lemma}) so that $\lim \tilde{f}(x_n) = \tilde{f}(x)$.
\end{proof}

We prove the generalization of theorem \ref{local_measurable_process_thm} as stated in remark \ref{measurable_process_remark}:
\begin{proof} (of generalization of \underline{theorem \ref{local_measurable_process_thm}}) We note that lemma \ref{uniform_lusin_measurable_lemma} is applicable. That is, fixing a probability measure $\mu$ on $S$ equivalent to Lebesgue measure, we get for each $\eps>0$ a measurable $K_{\eps}\subset S$ with $\mes(S\setminus K_{\eps})<\eps$ such that for all $t\in T$, $g(\cdot,t)|_{K_\eps}$ is continuous. We will assume in this proof, without loss of generality that every point of $K_\eps$ is a Lebesgue point of it (perhaps removing a set of measure zero from $K_\eps$).

Assume that $f$ is not a.e. equal to a function continuous on $\ld{S}+\od{T}$, then by lemma \ref{ae_measurable_continuity_criterion} there is a point $x_0\in \ld{S}+\od{T}$ and sequences $x_n^1\to x_0$ and $x_n^2\to x_0$ of \emph{Lebesgue points} of $f$ with different $f$ limits or with infinite $f$ limits. We now employ the same reasoning as in the proof of theorem \ref{local_measurable_thm} to get subsequences of $x_n^1$ and $x_n^2$ (which we continue to denote by $x_n^1$ and $x_n^2$) and $\delta>0$ such that 
\[T\cap\left[(x_0-K_\delta)\cap(\cap_n(x_n^1-K_\delta))\cap(\cap_n(x_n^2-K_\delta))\right]\ne \phi\]
Denote by $t_1\in T$ an element of this intersection. We claim that $g(\cdot,t_1)|_{K_\delta}$ is discontinuous and thus obtain a contradiction.

To see this, we use that $g(\cdot,t_1)$ is a.e. equal to $f$ and that every point of $K_\delta$ is a Lebesgue point of it to choose two sequences $x_n'^1$ and $x_n'^2$ each element of which is a small perturbation of $x_n^1$ and $x_n^2$ respectively with:
\begin{itemize}
\item $g(x,t_1)=f(x+t_1)$ for $x$ in these sequences.
\item The $f$ limits of these sequences are the same as the $f$ limits of $x_n^1$ and $x_n^2$.
\item Still $t_1\in T\cap\left[(x_0-K_\delta)\cap(\cap_n(x_n'^1-K_\delta))\cap(\cap_n(x_n'^2-K_\delta))\right]$.
\end{itemize}
This gives the required contradiction.
\end{proof}

\end{subsection}
\begin{subsection}{Some counter examples}
Before continuing to prove the other theorems, in order to emphasize the necessity of the conditions in the previous theorems, we briefly indicate some situations in which $f$ is $(S,T)$ shift continuous but not continuous on $\R$:
\begin{enumerate}
\item $S=\R$ and $T\ne\phi$ is countable, then any function $f:\R\to\R$ is $(S,T)$ shift continuous since we can take the discrete topology on $T$. 

Also if $T$ is any proper measurable subgroup of $\R$, define $f:=1_T$ then for every $s\in\R$ we have either $f(s+\cdot)|_T\equiv 0$ or $f(s+\cdot)|_T\equiv 1$, but $f$ is not continuous on $\R$.

If $S$ is countable and $T=\R$, then any function $f:\R\to\R$ is $(S,T)$ shift continuous since we can take the topology generated by $\{f(s+\cdot)\}_{s\in S}$ as the second countable topology.

Also, if $S$ is any proper measurable subgroup of $\R$, then $f:=1_S$ is $(S,\R)$ shift continuous since for any $s\in S$, $f(s+\cdot)=f(\cdot)$ (and we can take the topology generated by $f$).

We remark that proper measurable subgroups of $\R$ (it is easy to see that such a subgroup must have measure 0) exist with any Hausdorff dimension between 0 and 1 (see \cite{EV1}. Also related are \cite{EM1} and \cite{BO1} where it is shown that any Borel proper \emph{subring} has Hausdorff dimension 0).
\item $\mes(S^c)=0$ and $T$ is residual or $S$ is residual and $\mes(T^c)=0$. For the first case, take $S\in\R$ of full measure and first category and $T=(-S)^c$. Define $f(x):=\frac{1}{x}$. Then for every $s\in S$ we have $f(s+\cdot)|_T=\frac{1}{s+\cdot}$, taking the usual metric on $\R$ we get that $f$ is $(S,T)$ shift continuous (since $0\notin(S+T)$) but $f$ is not continuous on $\R$. The second case is very similar.
\item We note that if $f$ is continuous and $h$ is any homomorphism of the line, then $f(h(\cdot))$ is $(\R,\R)$ shift continuous (the second countable topology for $f\circ h$ being the $h$-image of the second countable topology for $f$).

For a counter example with $f$ not measurable, $S=T=\R$. Take $h$ to be a non-measurable homomorphism (solution of the Cauchy equation $h(x+y)=h(x)+h(y)$, for existence of such, see \cite{AC1} p. 35-36 where a construction of G. Hamel (1905) using a Hamel basis is described) and take for example $h$ itself (since $f(x)=x$ is continuous).

A similar example is given by $f:=1_G$ where $G$ is any proper subgroup of $\R$ with a countable number of cosets (it is easy to see that such a $G$ is never measurable. We can construct one as the span of a Hamel basis minus one element), since for any $s\in\R$, $f(s+G)=1_H$ where $H$ is one of the cosets, we can define a metric on $\R$ by taking the standard metric and setting distance 1 between any two cosets. Then $f$ is $(\R,\R)$ shift continuous, but clearly not continuous.
\end{enumerate}
The example before the last shows us that any measurable homomorphism is automatically continuous. While this is not new (see for example, thm 9.10 p. 61 \cite{KE1} (Baire property), \cite{FS1} (universally or Haar measurable in Abelian Polish groups), \cite{ZS1} (measurable in Locally compact groups) and \cite{AK1}, \cite{AK2} (measurable in locally compact groups, allowing change on measure zero)), it is interesting to note that this is a special case of our theorems when the target space of the homomorphism is second countable. Note, however, that more is known for homomorphisms then what follows from our theorems, for example, it is (in some situations) enough to assume only that the homomorphism has a measurable majorant to deduce that it is continuous, this is not the case for shift continuous functions in general as the last example shows.

Finally, it is clear that besides the case of a homomorphism, any measurable function satisfying a functional equation of the type $f(x,y)=F(f(x),f(y),x,y)$ (with $F$ a given continuous function) is also continuous from our theorems. This is also not new (see \cite{AC1} for more details).
\end{subsection}
\begin{subsection}{Stochastic processes application}
We finish this section by proving theorem \ref{stationary_process_thm}. 

We need the following well known fact:
\begin{lemma} \label{continuous_stationary_process}
A measurable stationary process $(X_t)_{t\in\R}$ is continuous in probability.
\end{lemma}
Since we couldn't find a reference for the proof, we briefly indicate it here.
\begin{proof}
First, since $L_0(\Omega)$ (see remark \ref{evaluation_function_remark}) is separable (even Polish) the set of random variables $(X_t)_{t\in\R}$ is separable, hence we can assume without loss of generality that the process generates the $\sigma$ algebra of $\Omega$ since otherwise we can pass to the quotient probability space $\Omega'$ which is still standard and for which the process does generate the $\sigma$ algebra. 

Now, since the process is stationary, for each shift $s\in\R$ we have an automorphism $T_s:\Omega\to\Omega$ of the probability space such that for each $t\in\R$, $X_{t+s}(\omega)=X_{t}(T_s(\omega))$ except on a set of 0 measure (not depending on $t$) and also $T_{s_1+s_2}(\omega)=T_{s_1}(T_{s_2}(\omega))$ a.s.. Assuming without loss of generality that $\Omega=[0,1]$ and denoting by $L_0^A(\Omega,\Omega)$ the set of equivalence classes of automorphisms of $[0,1]$ we easily check that this is a Polish group with respect to composition and convergence in probability \footnote{And hence this topology is unique (prop. 12.24 in \cite{KE1}).}. We consider the mapping $\alpha:\R\to L_0^A(\Omega,\Omega)$ defined by $\alpha(s):=T_s$. This is a measurable homomorphism (since the process $(X_t)_{t\in\R}$ is measurable and $(X_t)_{t\in\R}$ generate the $\sigma$-algebra of $\Omega$) and hence it is continuous (for example, by our theorem \ref{general_local_measurable_thm} in section \ref{general_formulation_section}). It remains to verify that if $T_{s_n}\to T_s$ in probability then for any random variable $f$, $f(T_{s_n})\to f(T_s)$ in probability and hence the process $(X_t)_{t\in\R}$ is continuous in probability.
\end{proof}

Continuing, We will use the following notation for sections of a product set. If $A\subset S\times T$ we write:
\[A_{(s,\cdot)}:=\{t\in T\ |\ (s,t)\in A\}\]
and similarly define $A_{(\cdot,t)}$.

We need a lemma which extends Fubini's theorem for a situation involving non-measurable subsets of $\R$ (we will need this when the set $T$ in the formulation of theorem \ref{stationary_process_thm} is not measurable), those readers not interested in non-measurable subsets can skip this lemma, assume that the set $T$ is measurable and use the standard Fubini theorem.

We denote the (Lebesgue) $\sigma$-algebra of $\R$ by $\cb$ and for any set $T\subset\R$ we define the induced $\sigma$-algebra by:
\[\cb_T:=\{B\cap T\ |\ B\in\cb\}\]
\begin{lemma} (extension of Fubini's theorem) \label{Fubini_extension}
Let $S,T\subset\R$ be arbitrary sets. Assume that $A\subset S\times T$ is measurable (in $\cb_S\otimes\cb_T$) and that there exists $T_0\subset T$ with $\mes(T\setminus T_0)=0$ such that for every $t\in T_0$ the section $A_{(\cdot,t)}$ has zero measure then there exists $S_0\subset S$ with $\mes(S\setminus S_0)=0$ such that for every $s\in S_0$ the section $A_{(s,\cdot)}$ has zero measure.
\end{lemma}
\begin{proof}
Let $T^0$ and $S^0$ be measurable covers of $T$ and $S$ respectively (cf. thm 3.3.1 in \cite{DU1}), that is, $T^0$ and $S^0$ are measurable, $T\subset T^0$, $S\subset S^0$ and $\mes^*(T)=\mes(T^0)$, $\mes^*(S)=\mes(S^0)$. Since $A$ is measurable in $\cb_S\otimes\cb_T$ there exists $A^1\subset\R^2$ measurable in $\R^2$ such that $A=A^1\cap(T_0\times S_0)$. Define $A^0:=A^1\cap(T^0\times S^0)$, we claim that $\mes_2(A^0)=0$. This follows by Fubini's theorem for $A^0$ since otherwise there would exist a \emph{measurable} set $T^{pos}\subset T^0$ with $\mes(T^{pos})>0$ such that for every $t\in T^{pos}$, $\mes(A^0_{(\cdot,t)})>0$. But from the assumptions, $T^{pos}\subset (T^0\setminus T_0)$ and hence $\mes_*(T^{pos})=0$, a contradiction. 

It now follows using Fubini's theorem for $A^0$ again that for a.e. $s\in S^0$, $\mes(A^0_{(s,\cdot)})=0$, hence in particular, the set of $s$ in $S_0$ which do not satisfy $\mes(A_{(s,\cdot)})=0$ is of measure 0.
\end{proof}
\begin{remark}
Note that it was not essential in the lemma that $T$ and $S$ were subsets of $\R$, the same proof gives the same result when $T\subset X_1$ and $S\subset X_2$ where $X_1$ and $X_2$ are arbitrary $\sigma$-finite measure spaces.
\end{remark}
\begin{proof} (of \underline{theorem \ref{stationary_process_thm}})
We denote by $f(\omega,t)$ a modification of $(X_t)_{t\in\R}$ measurable in both arguments jointly (such a modification exists since $(X_t)_t$ is measurable). Wishing to use theorem \ref{local_measurable_process_thm} we define a new process $(Y_t)_{t\in\R}$ on the probability space $(\Omega\times\R,\Prob\times\mu)$ where $\mu$ is any probability measure on $\R$ equivalent to Lebesgue measure by $Y_t(\omega,s):=f(\omega,s+t)$\footnote{One can think of two random draws, first draw a starting point in $\R$ according to $\mu$ and then draw a sample path according to $\Prob$.}. It is easy to see from the stationarity of $(X_t)_t$ that $(X_t)_t$ and $(Y_t)_t$ are identically distributed (that is, the joint distribution of the process in any finite number of times is identical). It follows from theorem 1(c) in \cite{TS1} that $(Y_t)|_T$ also has a natural modification which we denote by $g((\omega,s),t)$ (so that $g:(\Omega\times\R)\times T\to\R$). We note that $g$ is measurable in all variables jointly by theorem 2(d) in \cite{TS1} (with respect to the product sigma-algebra where $T$ is equipped with the induced $\sigma$-algebra $\cb|_T:=\{B\cap T\ |\ B\in\cb\}$. Note that the process $(X_t)|_T$ is still measurable with respect to this $\sigma$-algebra).

We note that since $f$ is measurable with respect to the product $\sigma$-algebra of $\Omega$ and $\R$, every section of it is measurable, that is, $f(\omega,\cdot)$ is measurable for every $\omega\in\Omega$. We define for each $\omega\in\Omega$ a "sub-process" of $(Y_t)_t$ on the probability space $(\R,\mu)$ and the time set $\R$ by $Z^\omega_t(s):=f(\omega, s+t)$. We will want to apply theorem \ref{local_measurable_process_thm} to it, to do that, we need to show that it has a natural modification.

We proceed to define two subsets of $\Omega$ with full measure:
\begin{itemize}
\item
First, since $g$ is a natural modification for $(Y_t)|_T$ there exists $A\subset\Omega\times\R$ with $(\Prob\times\mu)(A)=1$ and a second countable topology $\ct$ on $T$ such that for every $(\omega,s)\in A$, the function $g((\omega,s),\cdot)$ is $\ct$ continuous. From Fubini's theorem applied to $A$, there exists $\Omega_1\subset\Omega$ with $\Prob(\Omega_1)=1$ such that for every $\omega\in\Omega_1$, for $\mu$-a.e. $s\in\R$, the function $g((\omega,s),\cdot)$ is $\ct$ continuous.
\item
Second, we define 
\[C:=\{f(\omega,s+t)\ne g((\omega,s),t)\ |\ ((\omega,s),t)\in(\Omega\times\R)\times T\}\]
Then $C$ is measurable in $(\Omega\times\R)\times T$ (since $f$ and $g$ are) and for every $t\in T$, the section $C_{(\cdot,\cdot),t)}:=\{(\omega,s)\ |\ ((\omega,s),t)\in C\}$ has zero $(\Prob\times\mu)$ measure since $g$ is a modification of $(Y_t)_t$. Now using the extension of Fubini's theorem (lemma \ref{Fubini_extension}, or rather, its generalization in the remark after the lemma) we get that for a.e. $(\omega,s)\in\Omega\times\R$ the section $C_{((\omega,s),\cdot)}$ is null, again using Fubini's theorem there exists $\Omega_2\in\Omega$ with $\Prob(\Omega_2)=1$ such that for every $\omega\in\Omega_2$ for a.e. $s\in\R$, the section $C_{((\omega,s),\cdot)}$ is null and using the extension of Fubini's theorem again\footnote{Note that since $C$ is measurable in the product $\sigma$-algebra, for $\Prob$-a.e. $\omega\in\Omega$, the section $C_{((\omega,\cdot),\cdot)}$ is measurable in $\R\times T$.} we get that for every $\omega\in\Omega_2$ there exists $T_\omega$ with $\mes(T\setminus T_\omega)=0$ such that for every $t\in T_\omega$ the section $C_{((\omega,\cdot),t)}$ is $\mu$-null.
\end{itemize}
Define $\Omega^0:=\Omega_1\cap\Omega_2$ it follows that for $\omega\in\Omega^0$, the function $g((\omega,\cdot),\cdot)$ is a natural modification for $(Z^\omega_t)|_{T_\omega}$. Therefore, by theorem \ref{local_measurable_process_thm} for each such $\omega$, the function $f(\omega,\cdot)$ is a.e. equivalent to a continuous function which we denote by $h(\omega,\cdot)$.

We claim that $h(\cdot, \cdot)$ is a sample continuous modification of $(X_t)_{t\in\R}$. First, $h$ is measurable in both arguments jointly since it can be obtained from $f$ by using an evaluation function such as that described in remark \ref{evaluation_function_remark} with the integration applied on the $\R$ axis. Second, from Fubini's theorem we now get that for a.e. $t\in\R$, $h(\cdot,t)=f(\cdot,t)$ a.s., but the process $(X_t)_{t\in\R}$ is continuous in probability (lemma \ref{continuous_stationary_process}), hence for \emph{every} $t\in\R$, $h(\cdot,t)=f(\cdot,t)$ a.s. and we are done.
\end{proof}
\end{subsection}
\end{section}

\begin{section}{Generalization to Polish groups}\label{general_formulation_section}
In this section we will give general formulations of the theorems of section \ref{real_formulation_section} for Polish groups other than $\R$ and indicate what changes (if any) are needed for the proofs.

We will generalize our theorems to functions $f:G\to X$ where $G$ is a Polish group and $X$ is a separable metrizable space. As before, we will deduce that if a function is shift-continuous then it is continuous. This is a good place to remark the following to the reader:

Since the definition of shift-continuity does not require that $G$ be Polish, one can well ask what happens when the group $G$ is not so. We have not concerned ourselves much with this question, the reader should note that only when $G$ is second countable does it follow that every continuous function is shift-continuous, for example, when $G$ is metrizable but not separable with a left-invariant metric $d$ it is easy to check that the function $d(e,\cdot)$ ($e$ is the neutral element of $G$) is continuous but not shift-continuous since its shifts generate the entire topology (which is not separable). Also worth noting is that if $G$ is not Polish, but there exists a Polish topology $\ct$ on $G$ with the same Borel field (in \cite{KE1}, p. 81 such a $G$ is called \emph{Polishable}) then it follows immediately from our theorems that a function is shift-continuous on $G$ if and only if it is $\ct$-continuous (the original topology of $G$ does not play any role and in particular we get again prop. 12.24 p. 81 of \cite{KE1} that the topology $\ct$ is unique). For example, the shift-continuous functions on the group $C(\R)$ (continuous functions on $\R$ with the uniform convergence topology) are exactly those functions continuous relative to the uniform convergence on compact sets topology, in the same way, the shift-continuous functions on $\Z^\N$ with the discrete topology are exactly those functions continuous relative to the product topology.

We first generalize lemmas \ref{uniform_lusin_category_lemma} and \ref{uniform_lusin_measurable_lemma} which will be our main tools. Following \cite{TS1}, we have the following topological lemma, the lemma is proved there only for the case $Y=\R$ (see below), so we give a proof here\footnote{Our formulation is also slightly more general.}:
\begin{lemma} \label{uniform_lusin_topology_lemma}
Let $X$ be a second countable space, $Y$ a separable metrizable space, $C(X,Y)$ the space of continuous mapping from $X$ to $Y$. Then there exists a topology $\tau$ on $C(X,Y)$ with the following properties:
\begin{enumerate}
\item $\tau$ is second countable.
\item For any $x\in X$ the "evaluation mapping at the point x" $\phi_x:C(X,Y)\to Y$, $\phi_x(f):=f(x)$ is continuous relative to $\tau$.
\item \label{small_sigma_algebra_property} The Borel $\sigma$-algebra of $\tau$ equals the $\sigma$-algebra generated by $\{\phi_x\}_{x\in X}$.
\end{enumerate}
\end{lemma}
\begin{proof}
Let $\{U_n\}$ be a basis for the topology of $X$, $\{V_n\}$ a basis for the topology of $Y$, We define a sub-basis for the topology $\tau$ by:
\[B_m^n:=\{f\in C(X,Y)\ |\ f(U_m)\subset V_n\}\]
We claim that $\tau$ is the required topology. It is obvious that $\tau$ is second countable. Fix $x\in X$, Let us show that $\phi_x$ is continuous in $\tau$. Take any $V_n$ from the above basis for $Y$, it is clear that for any $m$ with $x\in U_m$, $B_m^n\subset\phi_x^{-1}(V_n)$. On the other hand, if $f\in\phi_x^{-1}(V_n)$ then from the continuity of $f$, there is some $U_m$ with $x\in U_m$ such that $f\in B_m^n$.

Denote the $\sigma$-algebra generated by $\{\phi_x\}_{x\in X}$ by $B_\phi$. To show property \ref{small_sigma_algebra_property} we note that since $\{\phi_x\}_{x\in X}$ are continuous the Borel $\sigma$-algebra of $\tau$ contains $B_\phi$. To show the reverse inclusion it is enough (since $\tau$ is second countable) to show that for each $m,n\in\N$, $B_m^n\in B_\phi$, and actually it will be enough to show that $\closure(B_m^n)\in B_\phi$ ($\closure$ denotes the closure operator), this is because $\closure(B_m^n)=\{f\in C(X,Y)\ |\ f(U_m)\subset\closure(V_n)\}$  and since $Y$ is metrizable $V_n=\cup_{\closure(V_k)\subset V_n} \closure(V_k)$. This in turn is clear since $f\in\closure(B_m^n)$ if and only if $f(x)\in\closure(V_n)$ for each $x$ in a dense countable subset of $U_m$.
\end{proof}
We now have the generalizations of lemmas \ref{uniform_lusin_category_lemma} and \ref{uniform_lusin_measurable_lemma}:
\begin{lemma} (generalization of lemma \ref{uniform_lusin_category_lemma}) \label{general_uniform_lusin_category_lemma}
Let $S$ be a metric space, $(T,\ct)$ a second countable space and $X$ a separable metrizable space. Given $g:S\times T\to X$ such that for each $t\in T$, $g(\cdot, t)$ has the Baire property and for each $s\in S$, $g(s,\cdot)$ is continuous then there exists a set $S_1\subset S$ with $S\setminus S_1$ of first category such that for each $t\in T$, $g(\cdot,t)|_{S_1}$ is continuous.
\end{lemma}
\begin{lemma} (generalization of lemma \ref{uniform_lusin_measurable_lemma}) \label{general_uniform_lusin_measurable_lemma}
Let $S$ be a metric space, $(T,\ct)$ a second countable space, $X$ a separable metrizable space and $\mu$ a probability measure on $S$. Given $g:S\times T\to X$ such that for each $t\in T$, $g(\cdot,t)$ is $\mu$-measurable\footnote{Measurable relative to the completion of $\mu$.} and for each $s\in S$, $g(s,\cdot)$ is continuous then the functions $\{g(\cdot,t)\}_{t\in T}$ are \emph{uniformly Lusin measurable} in the sense that for each $\eps>0$ there exists a closed $F_{\eps}\subset S$ with $\mu(S\setminus F_{\eps})<\eps$ such that for all $t\in T$, $g(\cdot,t)|_{F_\eps}$ is continuous\footnote{If $S$ is Polish, $F_\eps$ can be taken to be compact.}.
\end{lemma}
The proofs are as in \cite{TS1} (for more specific settings) and are nearly identical, we show here the proof of the first lemma:
\begin{proof} (of lemma \ref{general_uniform_lusin_category_lemma}) We endow the space of continuous functions $C((T,\ct), X)$ with the topology $\tau$ of lemma \ref{uniform_lusin_topology_lemma}. We consider $\tilde{g}:S\to C((T,\ct), X)$ defined by $\tilde{g}(s):=g(s,\cdot)$ We note that from the assumptions and the properties of $\tau$, $\tilde{g}$ is Baire measurable, hence from the second countability of $\tau$, there exists a set $S_1\subset S$ with $S\setminus S_1$ of first category such that $\tilde{g}|_{S_1}$ is continuous, again from the properties of $\tau$ this implies that for any $t\in T$, $g(\cdot,t)|_{S_1}$ is continuous.
\end{proof}

In order to generalize theorem \ref{local_measurable_thm} to $f:G\to X$ where $G$ is a Polish group we will need some definition of density point in $G$\footnote{When $G=\R^n$ or $G=\T^n$, theorem \ref{local_measurable_thm} is true without changes.}. While some definition is possible in a locally compact group $G$ (see lemma \ref{density_point_in_functions_lemma} below), in order to not introduce new notations, we will give a somewhat weaker version of the theorem where density points are not required.

\begin{theorem}[Theorem \ref{local_measurable_thm} in locally compact Polish groups] \label{general_local_measurable_thm}
Let $G$ be a locally compact Polish group with (left) Haar measure $\mu$, $X$ a separable metrizable space and $f:G\to X$ a (Haar) measurable function. If $f$ is $(S,T)$ shift continuous for a measurable $S\subset G$ and a set $T\subset G$ with either $\mu(S^c)=0$ and $\mu^*(T)>0$ or $\mu(S)>0$ and $\mu_*(T^c)=0$ then $f$ is continuous (at all points of $G$)\footnote{$\mu_*$ denotes inner measure.}.

If, moreover, $G$ is compact (and $\mu(G)=1$) then it is enough to require $\mu(S)+\mu^*(T)>1$ to deduce that $f$ is continuous.
\end{theorem}
\begin{proof}
The proof of this theorem is the same as that of theorem \ref{local_measurable_thm}, using lemma \ref{general_uniform_lusin_measurable_lemma} instead of lemma \ref{uniform_lusin_measurable_lemma}. Note also the following:

For any (Haar) measurable set $A\subset G$, the function $g\to\mu((gA)\cap A)$ is continuous by theorem 20.17 (p. 296) in \cite{HR1}.

We do not need to choose $s_0$ and $t_0$ as in the proof of theorem \ref{local_measurable_thm}, we take any $x_0\in G$ and the set $K_\delta^{n_k}$ with large enough measure so that it will intersect $T$\footnote{Note that $\mu(A)=0$ iff $\mu(A^{-1})=0$ by thm 15.14 (p. 197) in \cite{HR1}.}.
\end{proof}
More general Polish groups $G$ can be handled when the function $f$ has the Baire property (see theorem \ref{general_local_category_thm} below) or when it is universally measurable and the group is \emph{Abelian}. Unfortunately, we do not know a corresponding theorem in the non-Abelian case. Since we do not have a canonical measure on a general Polish group, the following theorem requires $(G,G)$ shift continuity:

\begin{theorem}[Theorem \ref{local_measurable_thm} in Abelian Polish groups] \label{global_universally_measurable_thm}
Let $G$ be an Abelian Polish group\footnote{More generally, any completely metrizable Abelian group (with the same proof).}, $X$ a separable metrizable space and $f:G\to X$ a universally measurable function. If $f$ is $(G,G)$ shift continuous then $f$ is continuous (in the metric of $G$).
\end{theorem}
The proof of this theorem is sufficiently different from that of theorem \ref{local_measurable_thm} that we give it in full.

Wishing to apply the method of proof from theorem \ref{local_measurable_thm}, we require some probability measure on $G$ with a property similar to the continuity of $\mes((x+A)\cap A)$. We have the following lemma which is enough for our needs\footnote{We do not know a similar lemma in the non-Abelian case and hence cannot prove theorem \ref{global_universally_measurable_thm} in that case.}:
\begin{lemma}  \label{Abelian_Polish_measure_lemma}
Let $G$ be an Abelian Polish group and $(x_n)_n\subset G$ a sequence which tends to 0. Then there exists a measure $\mu$ on $G$ such that for any $\mu$-measurable $A\subset G$ with $\mu(A)>0$ and $\eps>0$ there exists a subsequence $(x_{n_k})_k$ with $\mu((x_{n_k}+A)\cap A)>(1-\eps)\mu(A)$ for all $k$.
\end{lemma}
\begin{proof}
We start by defining a measure on the Polish space $\N^\N$, we define random variables $(X_n)_{n=1}^\infty$ with values in $\N$ by:
\[\Prob(X_n=j)=\alpha_n\left(\frac{n-1}{n}\right)^{-j}\]
with $\alpha_n$ the suitable normalizing constant. We consider the product space $(\N^\N,m)$ where $m$ is the product measure of the distributions of $(X_n)_n$. Defining the unit vectors $(e_n)_n\subset \N^\N$ we see that for any $m$-measurable $B\subset\N^\N$ we have:
\[m(B+e_n)=\frac{n-1}{n}m(B)\]
since we can calculate the measure of $B+e_n$ using the Fubini theorem when viewing $\N^\N$ as the product of the $n$'th coordinate and the rest of the coordinates.

Returning to the group $G$, we define a mapping $T:\N^\N\to G$ by
\[T((c_n)_n)=\sum_{n=1}^\infty c_n x_n\]
Since we are only interested in subsequences of $(x_n)_n$, we can assume without loss of generality that $(x_n)_n$ tends to 0 so quickly that the mapping $T$ will be defined $m$-a.s. on $\N^\N$ (using completeness of $G$), hence $T$ is $m$-measurable. We define $\mu$ by $\mu:=m\circ T^{-1}$.

Now, let $A\subset G$ be a $\mu$-measurable set with $\mu(A)>0$, let $\eps>0$. We consider $B:=T^{-1}(A)$. Fix  $\delta>0$ (to be chosen later), we can approximate $B$ by a finite union of basic open sets in $\N^\N$ (sets of the type $\prod_{i=1}^k V_i \times \prod_{j=k+1}^\infty \N$ with $V_i$ open in $\N$), we take $\{U_i\}_{i=1}^k$ such basic open sets whose union $U:=\cup_{i=1}^k U_i)$ satisfies:
\[m(B\cap U)\ge (1-\delta)m(U) \qquad \rm{ and } \qquad m(U)\ge(1-\delta)m(B)\]
Denote by $K$ the first coordinate from which there are no constraints in any of the $U_i$. Using the above mentioned property of $m$ we get for any $n\ge K$:
\[m((B+e_n)\cap U)=m((B+e_n)\cap (U+e_n))=\frac{n-1}{n}m(B\cap U)\]
And so we can take $N\ge K$ so large that for any $n\ge N$:
\[m((B+e_n)\cap U)\ge(1-2\delta)m(U)\]
And so:
\[m((B+e_n)\cap B)\ge (1-3\delta)m(U)\ge (1-3\delta)(1-\delta) m(B)\]
And by choosing the connection between $\delta$ and $\eps$ properly we get:
\[m((B+e_n)\cap B)\ge (1-\eps)m(B)\]
Returning now to the set $A$ in $G$, we finish by noting that\footnote{Here is the only place we use the assumption that $G$ is Abelian.}:
\begin{multline*}
\mu((A+x_n)\cap A)=m(T^{-1}((A+x_n)\cap A))=m(T^{-1}(A+x_n)\cap T^{-1}(A))\\
\ge m((B+e_n)\cap B) \ge (1-\eps) m(B)=(1-\eps) m(A)
\end{multline*}
As we wanted to prove.
\end{proof}
\begin{proof} (of \underline{theorem \ref{global_universally_measurable_thm}}) The proof is very similar to that of theorem \ref{local_measurable_thm}, we assume in order to get a contradiction that there exists a point $x_0\in G$ and a sequence $x_n\to x_0$ with $\lim f(x_n)\ne f(x_0)$, we then take the measure $\mu$ from lemma \ref{Abelian_Polish_measure_lemma} corresponding to the sequence $(x_n-x_0)_n$. Using uniform Lusin measurability (lemma \ref{general_uniform_lusin_measurable_lemma}) it follows that for every $\eps>0$ there exists $K_\eps\subset G$ with $\mu(G\setminus K_\eps)<\eps$ such that for all $t\in T$, $f(\cdot+t)|_{K_\eps}$ is continuous.

Now, fix any $\eps>0$. Using inductively the property of $\mu$ as given by lemma \ref{Abelian_Polish_measure_lemma} we find a subsequence $(x_{n_k})_k$ such that:
\[\mu((x_0-K_\eps)\cap(\cap_k(x_{n_k}-K_\eps)))>0\]
This is a contradiction, since for any $t\in ((x_0-K_\eps)\cap(\cap_k(x_{n_k}-K_\eps)))$, $f(\cdot+t)|_{K_\eps}$ is discontinuous.
\end{proof}

Switching to Baire category settings, the following generalizes theorem \ref{local_category_thm}:
\begin{theorem}[Theorem \ref{local_category_thm} in Polish groups] \label{general_local_category_thm}
Let $G$ be a Polish group\footnote{\label{Baire_space_footnote}More generally, any completely metrizable group or even any metric group which is a Baire space (p. 41 in \cite{MC1}), that is, for which each open set is of second category.}, $X$ a separable metrizable space and $f:G\to X$ a function with the Baire property. If $f$ is $(S,T)$ shift continuous for $S\subset G$ with the Baire property and a set $T\subset G$ then $f$ is continuous at all points of $\bd{S}+\sc{T}$.\footnote{Which is an open set since $\bd{S}$ is open.}
\end{theorem}
\begin{proof}
The proof is the same as that of theorem \ref{local_category_thm} using lemma \ref{general_uniform_lusin_category_lemma} instead of lemma \ref{uniform_lusin_category_lemma}.
\end{proof}
Note that remark \ref{large_category_remark} applies to this theorem also.

Moving on to theorems about stochastic fields we have the following generalization of theorem \ref{local_measurable_process_thm} to locally compact Polish groups:
\begin{theorem}[Theorem \ref{local_measurable_process_thm} in locally compact Polish groups] \label{general_local_measurable_process_thm}
Let $G$ be a locally compact Polish group with (left) Haar measure $\mu$, $X$ a separable metrizable space, $f:G\to X$ a (Haar) measurable function and $\Prob$ a probability measure on $G$ equivalent to $\mu$. Define a stochastic field $(X_t)_{t\in G}$ on $(G,\Prob)$ by $X_t(\omega):=f(\omega+t)$. Then there is an alternative:

Either $(X_t)_{t\in G}$ is sample continuous or for every non-null set $T\subset G$, $(X_t)|_T$ does not have a natural modification.
\end{theorem}
Note that remark \ref{sample_continuity_remark} is true here also.  To prove this theorem we need some notion similar to Lebesgue density points of real functions. Using this notion it is also possible to prove a local version of this theorem as in remark \ref{measurable_process_remark}\footnote{We can also define an analogue of outer density point in a similar manner.}, we will not formulate this generalization here.

A notion of density points adequate for our needs was suggested to us by Tsirelson \cite{TS2}. We first deal with sets $A\subset G$:
\begin{lemma} \label{density_point_in_sets_lemma}
Let $G$ be a locally compact metric group with (left) Haar measure $\mu$, $A\subset G$ a (Haar) measurable set, then for any sequence $\eps_n\to 0$ there is a subsequence $(\eps_{n_k})_k$ such that for a.e. $x\in A$:
\[\lim_{k\to\infty} \frac{\mu(A\cap B(x,\eps_{n_k}))}{\mu(B(x,\eps_{n_k}))}\to 1\]
where $B(x,\delta)$ is the (open) ball around $x$ of radius $\delta$ in some left invariant metric compatible with the topology.
\end{lemma}
For given $(\eps_n)_n$, we call a point where convergence in the above lemma occurs a \emph{density point of $A$ relative to $(\eps_{n_k})_k$}. Note that from the Lebesgue density theorem on $\R$ a.e. point of a measurable $A\subset\R$ is a density point of $A$ relative to any sequence.
\begin{proof}
We will assume the metric of the group to be left invariant (using thm 8.3 p. 70 in \cite{HR1}. The new metric might not be complete even if the topology admits a complete metric). 

It is enough to prove the lemma when $\mu(1_A)<\infty$, since $\mu$ is $\sigma$-finite, say $G\setminus N=\cup_n F_n$ with $N$ a $\mu$-null set and $F_n$ of finite $\mu$-measure.  We can apply the lemma to each $A\cap F_n$ each time refining the sequence $(\eps_n)_n$ more and more, and finally we can take the diagonal sequence of epsilons to get a sequence $(\eps_{n_k})_k$ good for a.e. $x\in A$.

Now, since $\mu(A)<\infty$, the indicator function $1_A(\cdot)$ is in $L_1(G)$, hence from theorem 20.15 (p. 293) in \cite{HR1} for each $\eps>0$ there is a $\delta>0$ such that:
\[\lVert 1_A * \left(\frac{1}{\mu(B(0,\delta_1))}1_{B(0,\delta_1)}\right) - 1_A\rVert_1<\eps \quad \forall \delta>\delta_1>0\]
It follows that we can choose a subsequence $(\eps_{n_k})_k$ such that 
\[f_k:=1_A * \left(\frac{1}{\mu(B(0,\eps_{n_k}))}1_{B(0,\eps_{n_k})}\right)\]
converges to $1_A$ a.s. (we first choose a subsequence to converge in $L_1(G)$ and then another subsequence to converge a.s.).  This finishes the proof since $f_k(x)=\frac{\mu(A\cap B(x,\eps_{n_k}))}{\mu(B(0,\eps_{n_k}))}$ and using the left invariance of $\mu$ and the metric.
\end{proof}
We immediately get a similar notion for functions:
\begin{lemma} \label{density_point_in_functions_lemma}
Let $G$ be a locally compact metric group with (left) Haar measure $\mu$, $X$ a second countable space and $f:G\to X$ a (Haar) measurable function. Then for any sequence $\eps_n\to 0$ there is a subsequence $(\eps_{n_k})_k$ such that for a.e. $x\in G$ and every open $U\subset X$ with $f(x)\in U$:
\[\lim_{k\to\infty} \frac{\mu(f^{-1}(U)\cap B(x,\eps_{n_k}))}{\mu(B(x,\eps_{n_k}))}\to 1\]
where $B(x,\delta)$ is as in lemma \ref{density_point_in_sets_lemma}.
\end{lemma}
For given $(\eps_n)_n$, we call a point where convergence in the above lemma occurs a \emph{density point of $f$ relative to $(\eps_{n_k})_k$}. Again, note that from the Lebesgue density theorem on $\R$ a.e. point of $\R$ is a density point of a measurable $f$ relative to any sequence.
\begin{proof}
We simply apply lemma \ref{density_point_in_sets_lemma} for each $U_n\subset X$ in a countable basis of the topology, each time refining the sequence $(\eps_n)_n$ more and more, and finally we take the diagonal sequence of epsilons $(\eps_{n_k})_k$. It is immediate that this subsequence satisfies the requirements of the lemma.
\end{proof}

Now the proof of theorem \ref{general_local_measurable_process_thm} is the same as that of theorem \ref{local_measurable_process_thm} but with Lebesgue density points of functions replaced by density points relative to some sequence tending to 0 (note that we can have the same sequence for even a countable number of functions, although in the proof we only need it for 2 functions: $f$ and $1_{K_\eps}$) and with lemma \ref{general_uniform_lusin_measurable_lemma} instead of lemma \ref{uniform_lusin_measurable_lemma}. Also note that $x_0$ in the proof should be just any $x_0\in G$.

The next theorem generalizes theorem \ref{stationary_process_thm}:
\begin{theorem}[Theorem \ref{stationary_process_thm} in locally compact Polish groups] \label{stationary_process_locally_compact_thm}
Let $G$ be a locally compact Polish group with (left) Haar measure $\mu$, $X$ a separable metrizable space. Let $(X_t)_{t\in G}$ be a (Haar) measurable \emph{stationary} stochastic field with values in $X$ (i.e., for each $t\in G$, $X_t:\Omega\to X$). Then the following alternative holds:

Either $(X_t)_{t\in G}$ is sample continuous or for every non-null set $T\subset G$, $(X_t)|_T$ does not have a natural modification.
\end{theorem}
Note that remark \ref{second_sample_continuity_remark} is true here also. The proof is the same as that of theorem \ref{stationary_process_thm}, for an evaluation function (for the end of the proof) one can use, for example, $\alpha:G\times L_0^C(G)\to X$ defined by $\alpha(x,f):=\tilde{f}(x)$ (where $\tilde{f}$ is as defined in lemma \ref{ae_measurable_continuity_criterion}, density points are as in lemma \ref{density_point_in_functions_lemma} and we denoted by $L_0^C(G)$ the equivalence classes (mod 0) of continuous functions from $G$ to $X$).

We continue with the generalization of theorem \ref{local_category_process_thm}:
\begin{theorem}[Theorem \ref{local_category_process_thm} in Polish groups] \label{local_category_process_polish_groups_thm}
Let $G$ be a Polish group\footnote{More generally, any completely metrizable group or even any metric group which is a Baire space (see footnote \ref{Baire_space_footnote}).}, $X$ a separable metrizable space and $f:G\to X$ a function with the Baire property. Given a set $T\subset G$ of second category and a function $g:G^2\to X$ such that for each $t\in T$ $g(\cdot,t)$ is equal to $f(\cdot+t)$ on a residual set, then if  there exists a second countable topology $\ct$ on $T$ such that for a residual set of $s\in G$, $g(s,\cdot)|_T$ is $\ct$-continuous, then $f$ is equal to a continuous function $h$ on a residual set and for a residual set of $s\in G$, $g(s,\cdot)|_T=h(s+\cdot)|_T$.
\end{theorem}
Note that remark \ref{local_category_process_remark} is true here also. The proof is the same as that of theorem \ref{local_category_process_thm} using lemma \ref{general_uniform_lusin_category_lemma} instead of lemma \ref{uniform_lusin_category_lemma}.

We end by giving a topological corollary of our work which appears to be new:
\begin{corollary}
Let $G$ be a Polish group and denote its topology by $\ct$. If $\ct'$ is a separable metrizable topology on $G$ such that for each $g\in G$ the group operation $g\cdot h$ is continuous in $h\in G$ with respect to $\ct'$ and the sets in $\ct'$ are either Baire measurable, Haar measurable (in locally compact groups) or universally measurable (in Abelian Polish groups) then $\ct'\subset\ct$.
\end{corollary}
\begin{proof}
Let $id:(G,\ct)\to(G,\ct')$ be the identity map on $G$, then $id$ is measurable (in the same sense that the sets in $\ct'$ are) and $(G,G)$ shift continuous (since for each $g\in G$, $g\cdot h$ is continuous in $h\in G$ with respect to $\ct'$), hence from our theorems (thm \ref{general_local_measurable_thm}, \ref{global_universally_measurable_thm} or \ref{general_local_category_thm}) $id$ is continuous, that is, $\ct'\subset\ct$.
\end{proof}
\end{section}
\begin{section}{Open problems}
In this section we formulate some questions related to our work which are still open:
\begin{enumerate}
\item \label{stationary_process_polish_group_problem} The main open question (from our point of view) is whether theorem \ref{stationary_process_thm} can be generalized to arbitrary Polish groups, i.e., is any measurable stationary stochastic field on a Polish group $G$ sample continuous if and only if it has a natural modification? All we know is what is written in theorem \ref{stationary_process_locally_compact_thm}, that is, when the group $G$ is locally compact. It would be interesting to try and prove this (or find a counter example) when $G$ is a Hilbert or a Banach space.
\item Can theorem \ref{global_universally_measurable_thm} be generalized to non-Abelian groups? That is, is any universally measurable function on a (non-Abelian) Polish group $G$ continuous if and only if it is $(G,G)$ shift continuous? All we require for our methods to work in such a group is a construction of a measure on the non-Abelian group with the properties of the measure in lemma \ref{Abelian_Polish_measure_lemma}. In fact, it would be enough that given $x_n\to 0$ in $G$, there exists a measure on $G$ such that for any set of positive measure $A$ corresponds a subsequence $(x_{n_k})_k$ with:
\[(-A)\cap\left(\cap_k(-A+x_{n_k})\right)\ne\phi\]
Unfortunately, we don't know that such a measure exists.
\item In connection with problem \ref{stationary_process_polish_group_problem} above, our proof of theorem \ref{stationary_process_thm} was by reducing the problem to that of theorem \ref{general_local_measurable_process_thm}, i.e., we showed that a.e. sample path of the stationary process was (up to a change of measure zero of each shift) shift continuous and deduced from that that the path is actually continuous.

Can theorem \ref{general_local_measurable_process_thm} be generalized to arbitrary Polish groups? In order to state the theorem in a general Polish group $G$ we need the notion of a \emph{null set} in $G$, we know of 4 such notions in general Polish groups: A set of first category, Christensen's Haar zero sets (in Abelian Polish groups, see \cite{CH1}, \cite{HSY1} and \cite{HSY2}), Aronszajn's null sets (in Banach spaces. Also called Gaussian null sets, see \cite{AR1}, \cite{PH1} and \cite{CS1}) and the relatively new Lindenstrauss and Preiss' $\Gamma$-null sets (also in Banach spaces, see \cite{LP1}). Given any such notion, one possible formulation of the question is:

Let $G$ be a Polish group, $f:G\to \R$ a Borel measurable function. Given $g:G^2\to\R$ such that for each $t\in G$, $g(\cdot,t)$ is a.e. equal to $f(\cdot+t)$. Does existence of a second countable topology $\ct$ on $G$ such that for each $s\in G$, $g(s,\cdot)$ is $\ct$-continuous imply that $f$ is a.e. equal to a continuous function?

Where we mean a.e. equal in the sense of the null sets. For first category sets, this is proved in theorem \ref{local_category_process_polish_groups_thm}. For the other notions, we don't know. Our proof requires a form of Lusin's theorem (that the restriction of $f$ to a large set is continuous) and some notion of density point of the function $f$ agreeing with the notion of null set. 

It might be useful to first think of the following (related) simpler question which only uses the original metric of the group:

Let $G$ be a Polish group, $f:G\to \R$ a Borel measurable function and $S\subset G$ a non-null Borel set. Assume that for every $t\in G$ the function $f(\cdot+t)|_S$ can be changed on a null set so that this restriction is continuous, does it follow that $f$ is a.e. equal to a continuous function?

Here, again, some notion of density point of $f$ agreeing with the notion of null set would be useful\footnote{We would like that a.e. point of a non-null Borel set $A$ has the property that each neighborhood of it has a non-null intersection with $A$ and also that if $x$ is a such a point for both $A$ and $B$, then it is such a point for $A\cap B$ also.}. 

Solution of these problems may help solve problem \ref{stationary_process_polish_group_problem} above. Note however, that in our proof of theorem \ref{stationary_process_locally_compact_thm} we also made strong use of Fubini's theorem. Unfortunately, In \cite{CH1} a counter example to Fubini's theorem is given for the class of Christensen's Haar zero sets, the example also works for Aronszajn's null sets (we don't know if Fubini's theorem is valid for $\Gamma$-null sets), so our specific method of proof cannot work using these two notions. 
\end{enumerate}
\end{section}

\begin{section}{Acknowledgments}
We wish to thank Prof. B. Weiss for some useful insights concerning the measure in lemma \ref{Abelian_Polish_measure_lemma}. We also wish to thank Mr.\!\! Nir Lev for useful discussions of related problems.
\end{section}

\begin{section}{Appendix: Olevskii's theorem}
In this section we give a different proof of theorem \ref{local_measurable_process_thm} using a theorem of A. Olevskii \cite{OL1}. The theorem and proof are very beautiful in their own right. They were not given in the text of the article since a simpler proof was found which is slightly more general (Olevskii's theorem does not allow the set $T$ to be non-measurable).
We begin by phrasing and proving Olevskii's theorem and then we comment on how it is equivalent to our theorems (for the cases of a measurable set $T$).
\begin{theorem} (A. Olevskii) \label{Olevskii_theorem}
Let $f:\R\to\R$ be a measurable function and $A\subset\R$ measurable with positive measure. If $\{f(\cdot+t)\}_{t\in\R}$ is separable in $L_\infty(A)$ then $f$ is a.e. equal to a continuous function.
\end{theorem}
\begin{proof}
We will need the following remarks:

For a function $g:\R\to\R$ we denote by $g_t$ the \emph{shifted} function $g_t(s):=g(s+t)$.

If $g$ and $h$ are in $L_\infty(A)$ and $x$ is a Lebesgue density point of $g,h$ and $A$ then 
\[d_{L_\infty(A)}(g,h)\ge|g(x)-h(x)|\] 
and we note in particular that both $g(x)$ and $h(x)$ are well defined and finite. It also follows that for every $\eps>0$, there exists a $\delta>0$ such that
\[d_{L_\infty(A)}(g_{t_1},h_{t_2})\ge(1-\eps)|g(x)-h(x)| \quad\ \forall |t_1|,|t_2|<\delta\]

Starting the proof, we assume without loss of generality that all points of $A$ are Lebesgue density points. We also note that $f$ is locally bounded in $L_{\infty}$ sense (otherwise there is a point $x_0$ in which $f$ is not locally bounded, any rotation of $f$ which "puts" $x_0$ on a point of $A$ will not be in $L_{\infty}(A)$).

We assume that $f$ is not a.e. equal to a continuous function, then since $f$ is locally bounded it follows from lemma \ref{ae_measurable_continuity_criterion} that there exists $x_0$ and sequences $x^1_n$ and $x^2_n$ of Lebesgue density points of $f$ such that $x^1_n\to x_0$ and $x^2_n\to x_0$ but $f(x^1_n)$ and $f(x^2_n)$ tend to different limits. We assume without loss of generality that $x_0\in A$ (by shifting $f$) and that $f(x^1_n)\to1$ and $f(x^2_n)\to-1$.

The following lemma is the heart of the proof:
\begin{lemma} \label{Olevskii_lemma}
For every $\eps>0$, there exist shifts $t_1,t_2$ and $\delta>0$ such that $|t_1|,|t_2|<\eps$, $(x_0-t_1),(x_0-t_2)\in A$ and $d_{L_{\infty}(A)}(f_{t_1+r_1},f_{t_2+r_2})>1$ whenever $|r_1|,|r_2|<\delta$.
\end{lemma}
In words, for every $\eps>0$ we can find two shifts smaller than $\eps$ which (each) move the discontinuity point of $f$ to another Lebesgue density point of $A$ and for which the shifted functions and their small shifts are at distance larger than one from each other.
\begin{proof}
We choose $n$ large enough so that $|x^1_n-x_0|,|x^2_n-x_0|<c\eps$ and $f(x^1_n)>\frac{1}{2}$, $f(x^2_n)<-\frac{1}{2}$ for a constant $1>c>0$ to be determined later. We claim that (for $c$ small enough) there exist $t_1,t_2$ with $|t_1|,|t_2|<\eps$ which satisfy:
\begin{enumerate}
\item $x^1_n-t_1=x^2_n-t_2$ and $(x^1_n-t_1)\in A$.
\item $(x_0-t_1)\in A$ and $(x_0-t_2)\in A$.
\end{enumerate}
Noting that in particular $t_2=x^2_n-x^1_n+t_1$ and writing $A_1:=A-x_0$ we see that the first property is satisfied for any $t_1\in(-A_1+x_n^1-x_0)$ and the second property for any $t_1\in(-A_1\cap(-A_1+x_n^1-x_n^2))$. Now since $x_0$ is a Lebesgue density point of $A$,  we see that the probability of $t_1$ chosen uniformly in $(-c\eps,c\eps)$ to satisfy both properties tends to 1 as $c$ tends to 0.

The lemma now follows when the remarks at the beginning of the proof are applied to $g:=f_{t_1}, h:=f_{t_2}$ and $x:=x^1_n-t_1$.
\end{proof}
The lemma contradicts separability of $\{f_t\}_{t\in\R}$ in $L_{\infty}(A)$. To see this we build a "Cantor" set $C\subset\R$ (with a continuum) of shifts such that for every different $t_1,t_2\in C$ we have $d_{L_{\infty}(A)}(f_{t_1},f_{t_2})>1$, this set is constructed by repeated use of the lemma, first for $f$ with $\eps:=1$ to produce two intervals $I^1_1:=[t_1-\eps_1,t_1+\eps_1]$ and $I^1_2:=[t_2-\eps_1, t_2+\eps_1]$ with the property that $d_{L_{\infty}(A)}(f_{r_1},f_{r_2})>1$ for $r_1\in I^1_1, r_2\in I^1_2$. Then we apply the lemma to $f_{t_1}$ with $\eps:=\eps_1$ to produce two sub-intervals $I^2_1,I^2_2\subset I^1_1$ (we use that $(x_0-t_1)\in A$ in order to apply the lemma) and similarly to $f_{t_2}$ with the same $\eps$ to produce $I^2_3,I^2_4\subset I^1_2$. Continuing in this manner we get closed intervals $(I^n_j)_{j=1\ldots 2^n}^n$ with $I^n_j\subset I^{n-1}_{\lfloor (j+1)/2\rfloor}$ and with $d_{L_{\infty}(A)}(f_{r_1},f_{r_2})>1$ for $r_1\in I^n_j, r_2\in I^n_k, j\ne k$. We define $I^n:=\cup_j I^n_j$ and $C:=\cap I^n$. It is easy to see that $C$ has the desired properties.
\end{proof}
Examining the proof more carefully, we see that it implies the following local version:
\begin{corollary} \label{local_separability_thm}
Let $f:\R\to\R$ be a measurable function and $A,T\subset\R$ measurable sets with positive measure. If $\{f(\cdot+t)\}_{t\in T}$ is separable in $L_\infty(A)$, then $f$ is a.e. equal to a function $g$ continuous on $\ld{A}+\ld{T}$.
\end{corollary}
The equivalence to theorem \ref{local_measurable_process_thm} (and even its generalization in remark \ref{measurable_process_remark}) now follows from the results of Tsirelson (thm 1(e) in \cite{TS1}) which show (in our situation) that for $S,T\subset\R$, a process $(X_t)_{t\in T}$ defined on a (standard) probability space $(S,\Prob)$ by $X_t(\omega):=f(\omega+t)$ has a natural modification if and only if for any $\eps>0$ there exists $A\subset S$ with $\Prob(A)\ge 1-\eps$ such that $\{f(\cdot+t)\}_{t\in T}$ is separable in $L_\infty(A)$.

We remark also that we used one side of the equivalence in our proof when we established uniform Lusin measurability of $\{f(\cdot+t)\}_{t\in t}$ (lemma \ref{uniform_lusin_measurable_lemma}). The other side of the equivalence (from separability in $L_\infty$ to uniform Lusin measurability) is established by Lusin's theorem (see \cite{TS1} for more details).
\end{section}

\end{document}